\documentclass[12pt]{amsart}
\usepackage{amsmath,amsthm,amsfonts,amssymb}
\begin{document} 
\newcommand{\A}{{\mathbb A}}
\newcommand{\B}{{\mathbb B}}
\newcommand{\C}{{\mathbb C}}
\newcommand{\N}{{\mathbb N}}
\newcommand{\Q}{{\mathbb Q}}
\newcommand{\Z}{{\mathbb Z}}
\renewcommand{\P}{{\mathbb P}}
\newcommand{\R}{{\mathbb R}}
\newcommand{\rc}{\subset}
\newcommand{\rank}{\mathop{rank}}
\newcommand{\trace}{\mathop{tr}}
\newcommand{\dimc}{\mathop{dim}_{\C}}
\newcommand{\Lie}{\mathop{Lie}}
\newcommand{\Spec}{\mathop{Spec}}
\newcommand{\Auto}{\mathop{{\rm Aut}_{\mathcal O}}}
\newcommand{\alg}[1]{{\mathbf #1}}
\newtheorem*{definition}{Definition}
\newtheorem*{claim}{Claim}
\newtheorem{corollary}{Corollary}
\newtheorem*{Conjecture}{Conjecture}
\newtheorem*{SpecAss}{Special Assumptions}
\newtheorem{example}{Example}
\newtheorem*{remark}{Remark}
\newtheorem*{observation}{Observation}
\newtheorem*{question}{Question}
\newtheorem*{fact}{Fact}
\newtheorem*{remarks}{Remarks}
\newtheorem{lemma}{Lemma}
\newtheorem*{proposition}{Proposition}
\newtheorem{theorem}{Theorem}
\title{%
On tameness and growth conditions
}
\author {J\"org Winkelmann}
\begin{abstract}
We study discrete subsets of $\C^d$, relating ``tameness'' with
growth conditions.
\end{abstract}
%
%
\address{%
J\"org Winkelmann \\
Mathematisches Institut\\
Universit\"at Bayreuth\\
Universit\"atsstra\ss e 30\\
D-95447 Bayreuth\\
Germany\\
}
\email{jwinkel@member.ams.org\newline\indent{\itshape Webpage: }%
http://btm8x5.mat.uni-bayreuth.de/\~{ }winkelmann/
}
\maketitle

\section{Results}

A discrete subset $D$ in $\C^n$ ( $n\ge 2$) is called
``tame'' if there exists a holomorphic automorphism $\phi$
of $\C^n$ such that $\phi(D)=\Z\times\{0\}^{n-1}$
(see  \cite{RR}).
If there exists a linear projection $\pi$ of $\C^n$ onto some
$\C^k$ ($0<k<n$) for which the image $\pi(D)$ is discrete,
then $D$ is tame (\cite{RR}). 
If $D$ is a discrete subgroup (e.g.~a lattice)
of the additive group $(\C^n+)$, then $D$ must be tame
(\cite{B}, lemma 4.4 in combination with corollary~2.6).
On the other hand there do exist discrete
subsets which are not tame (see \cite{RR}, theorem~3.9).

Here we will investigate how ``tameness'' is related to growth
conditions for $D$.

Slow growth implies tameness as we well see.
On the other hand, rapid growth can not imply non-tameness,
since every discrete subset of $\C^{n-1}$ is tame regarded
as subset of $\C^n=\C^{n-1}\times\C$.

The key method is to show that sufficiently slow growth implies that
a generic linear projection will have discrete image for $D$.

The main result is:

\begin{theorem}\label{thm-main}
Let $n$ be a natural number and let
$v_k$ be a sequence of elements in $V=\C^n$.

Assume that
\[
\sum_k \frac{1}{||v_k||^{2n-2}}<\infty
\]

Then $D=\{v_k:k\in\N\}$ is {\em tame}, i.e., 
there exists a biholomorphic map
$\phi:\C^n\to\C^n$
such that
\[
\phi(D)=\Z\times\{0\}^{n-1}.
\]
\end{theorem}

This growth condition is fulfilled for discrete subgroups of rank
at most $2n-3$, implying the following well-known fact:
\begin{corollary}
Let $\Gamma$ be a discrete subgroup of $\Z$-rank at most $2n-3$
of the additive group $(\C^n,+)$.

Then $\Gamma$ is a tame discrete subset of $\C^n$.
\end{corollary}

While this is well-known (even with no condition on the $\Z$-rank of
$\Gamma$), our approach yields the additional information that these
discrete subsets remain tame after a small deformation:

\begin{corollary}
Let $\Gamma$ be a discrete subgroup of $\Z$-rank at most $2n-3$
of the additive group $(\C^n,+)$,
$0<\lambda<1$ and $K>0$. Let $D$ be a subset of $\C^n$
for which there exists a bijective map $\zeta:\Gamma\to D$
with
\[
||\zeta(v)-v|| \le \lambda||v||+K
\]
for all $v\in\Gamma$.

Then $D$ is a tame discrete subset of $\C^n$.
\end{corollary}

This confirms the idea that tame sets should be stable under
deformation.
Similarily one would hope that the category of non-tame sets
is also stable under deformation. Here, however, one has to be
careful not to be too optimistic, because in fact the following is true:

\begin{proposition}
For every non-tame discrete subset $D\subset\C^n$ ($n>1$)
there is a tame discrete subset $D'$ and a bijection $\alpha:D\to D'$
such that
\[
||\alpha(v)-v|| \le \frac{1}{\sqrt 2}||v|| \ \ \forall v\in D
\]
and
\[
||w-\alpha^{-1}(w)|| \le ||w|| \ \ \forall w\in D'.
\]
\end{proposition}

In particular, if $D$ is a tame discrete subset and
$\zeta:D\to\C^n$ is a bijective map with
$||\zeta(v)-v||\le||v||$ for all $v\in D$, it is
possible that $\zeta(D)$ is not tame.

Still, one might hope for a positive answer to the following question:
\begin{question}
Let $n\in\N$, $n\ge 2$, let $1>\lambda>0$, $K>0$, let $D$ be a tame
discrete subset of $\C^n$ and let $\zeta:D\to\C^n$ be a map
such that
\[
||\zeta(v)-v|| \le \lambda ||v|| +K
\]
for all $v\in D$.
Does this imply that $\zeta(D)$ is a {\em tame} discrete subset of
$\C^n$ ?
\end{question}

Technically, the following is the key point for the proof
of our main result (theorem~\ref{thm-main}):

\begin{theorem}\label{thm-tech}
Let $n>d>0$.
Let $V$ be a complex vector space of dimension $n$ and let
$v_k$ be a sequence of elements in $V$.

Assume that
\[
\sum_k \frac{1}{||v_k||^{2d}}<\infty
\]

Then there exists a complex linear map $\pi:V\to\C^d$ such that
the set of all $\pi(v_k)$ is discrete in $\C^d$.
\end{theorem}

In a similar way on can prove such a result for real vector spaces:

\begin{theorem}\label{thm-real}
Let $n>d>0$.
Let $V$ be a real vector space of dimension $n$ and let
$v_k$ be a sequence of elements in $V$.

Assume that
\[
\sum_k \frac{1}{||v_k||^{d}}<\infty
\]

Then there exists a real linear map $\pi:V\to\R^d$ such that
the set of all $\pi(v_k)$ is discrete in $\R^d$.
\end{theorem}


For the proof of the existence of a linear projection $\pi$
with $\pi(D)$ discrete we proceed by
regarding randomly chosen
linear projections and verifying that the image of $D$ under
a random projection has discrete image with probability $1$
if the above stated series converges.

\section{Proofs}

First we deduce an auxiliary lemma.
\begin{lemma}\label{l-sphere}
Let $k,m>0$, $n=k+m$ and let $S$  denote the unit sphere
in $\R^n=\R^k\oplus\R^m$.
Furthermore let
\[
M_\epsilon= 
\{ (v,w)\in \R^k\times\R^m: ||v||\le\epsilon, (v,w)\in S\}.
\]
Then there are constants $\delta>0$, $C_1>C_2>0$ such that for
all $\epsilon<\delta$ we have
\[
C_1\epsilon^{k}\ge \lambda(M_\epsilon)\ge C_2\epsilon^{k}
\]
where $\lambda$ denotes the rotationally invariant probability measure
on $S$.
\end{lemma}
\begin{proof}
For each $\epsilon\in ]0,1[$ there is a bijection
\[
\phi_\epsilon: B\times S'\to M_\epsilon
\]
where
\[
B=\{v\in\R^k:||v||\le 1\},\quad S'=\{w\in \R^m:||w||= 1\}
\]
and
\[
\phi_\epsilon(v,w)=\left(\epsilon v;
\sqrt{1-||\epsilon v||^2}w\right).
\]
The functional determinant for $\phi_\epsilon$
equals 
\[
\epsilon^k\left(\sqrt{1-||\epsilon v||^2}\right)^m.
\]
It follows that
\[
\epsilon^k\left(\sqrt{1-\epsilon^2}\right)^m 
volume(S'\times B)\le
volume(M_\epsilon)
\le
\epsilon^k volume(S'\times B),
\]
which in turn implies
\[
\lim_{\epsilon\to 0}
\epsilon^{-k}\frac{volume(M_\epsilon)}{volume(S'\times B)}=1.
\]
Hence the assertion.
\end{proof}

\begin{lemma}
Let $\Gamma$ be a discrete subgroup of $\Z$-rank $d$ in $V=\R^n$.

Then
\[
\sum_{\gamma\in\Gamma}||\gamma||^{-d-\epsilon}<\infty
\]
for all $\epsilon>0$.
\end{lemma}

\begin{proof}
Since all norms on a finite-dimensional vector space are
equivalent, there is no loss in generality if we assume that
the norm is the maximum norm and $\Gamma=\Z^d\times\{0\}^{n-d}$.
Then the assertion is an easy consequence of the fact that
$\sum_{n\in\N} n^{-s}<\infty$ if and only if $s>1$.
\end{proof}

Now we proceed with the proof of theorem~\ref{thm-tech}:

\begin{proof}
We fix a surjective  linear map $L:V\to W=\C^d$.
Let $K$ denote $U(n)$ (the group of unitary complex linear transformations
of $V$).
For each $g\in K$ we define a linear map $\pi_g:V\to W$
as follows:
\[
\pi_g:v\mapsto L(g\cdot v).
\]
For $k\in\N$ and $r\in\R^+$ define
\[
S_{k,r}=\{g\in K: ||\pi_g(v_k)||\le r\},
\]
\[
M_{N,r}=\{g\in K: \#\{k\in\N: g\in S_{k,r}\}\ge N\}
\]
and
\[
M_r=\cap_N M_{N,r}.
\]
Now for each $g\in K$ the set $\{\pi_g(v_k):k\in\N\}$ is discrete
unless there is a number $r>0$ such that infinitely many distinct
image points are contained in a ball of radius $r$.
By the definition of the sets $M_r$ it follows that 
$\{\pi_g(v_k):k\in\N\}$ is discrete unless
$g\in M=\cup M_r$.

Let us now assume that there is no linear map $L':V\to W$
with $L'(D)$ discrete. Then $K=M$. In particular $\mu(M)>0$,
where $\mu$ denotes the Haar measure on the compact topological
group $K$. Since the sets $M_r$
are increasing in $r$, we have
\[
M=\cup_{r\in\R^+}M_r=\cup_{r\in\N}M_r
\]
and may thus deduce that $\mu(M_r)>0$ for some number $r$.
Fix such a number $r>0$ and define $c=\mu(M_r)>0$.
Then $\mu(M_{N,r})\ge c$ for all $N$, since $M_r=\cap M_{N,r}$.
However, for fixed $N$ and $r$ we have
\[
N\mu(M_{N,r})\le \sum_k\mu(S_{k,r}).
\]
Hence
\[
\sum_{k\in\N}\mu(S_{k,r})\ge N\mu(M_{N,r})\ge Nc
\]
for all $N\in\N$. Since $c>0$, it follows
that
$\sum_k\mu(S_{k,r})=+\infty$.

Let us now embedd $\C^d$ into $\C^n$ as the orthogonal complement
of $\ker L$. In this way we may assume that $L$ is simply the map
which projects a vector  onto its first $d$
coordinates, i.e., 
\[
L(w_1,\ldots,w_n) = (w_1,\ldots,w_d;0,\ldots,0).
\]
Now $g\in S_{k,r}$ is equivalent to the condition that $g(v_k)$
is a real multiple of an element in $M_\epsilon$ where $M_\epsilon$
is defined as in lemma~\ref{l-sphere} with $\epsilon=r/||v_k||$.
Using lemma~\ref{l-sphere} we may deduce that
$\sum_k\mu(S_{k,r})$ converges if and only if
$\sum_k ||v_k||^{-2d}$ converges.
\end{proof}

\begin{proof}[Proof of theorem~1]
The growth condition allows us to employ theorem~\ref{thm-tech}
in order to deduce that there is a linear projection onto a 
space of complex dimension $d-1$ which maps $D$ onto a discrete image.
By the results of Rosay and Rudin it follows that $D$ is tame.
\end{proof}

\begin{proof}[Proof of the proposition]
We fix a decomposition $\C^n=\C\times\C^{n-1}$ and write $D$ as the
union of all $(a_k,b_k)\in\C\times\C^{n-1}$ ($k\in\N$).
We define 
\[
\alpha(a_k,b_k)=
\begin{cases}
(a_k,0) & \text{ if }||a_k||>||b_k|| \\
(0,b_k) & \text{ if }||a_k||\le||b_k|| \\
\end{cases}
\]
Then $D'=\alpha(D)$ is tame because each of the projections to one of
the two factors $\C$ and $\C^{n-1}$ maps $D'$ onto a discrete subset.

The other assertions follow from the triangle inequality.
\end{proof}

The proof of thm.~\ref{thm-real} works in the same way as the
proof of thm.~\ref{thm-tech}, simply using the group of all
orthogonal transformations instead of the group of unitary
transformations.


\begin{thebibliography}{Bla}
\bibitem{B} Buzzard, G.:
Tame sets, dominating maps and complex tori.
\sl Trans. A.M.S. \bf 355\rm, 2557-2568 (2002)

\bibitem{BF} Buzzard, G.; Forstneric, F.:
An interpolation theorem for holomorphic automorphisms of $\C^n$. 
\sl J.~Geom. Anal. \bf 10\rm, 101-108 (2000)

\bibitem{RR} Rosay, J.P.; Rudin, W.: Holomorphic maps from $\C^n$ to
$\C^n$. \sl Trans. A.M.S. \bf 310\rm, 47--86 (1988)
\end{thebibliography}
\end{document}